% Logic Eprints
%Submitted 0426 Tue Jun 20, 1995 by: gitik@math.tau.ac.il (gitik moti )
%logic/gitik/club-reg.tex
%

%THIS IS A REWRITE OF DIANA's PAPER WHICH I
%STARTED ON DECEMBER 29th, 1993.

%THIS IS AGAIN A REWRITE WHICH I STARTED ON
%OCTOBER 31,1994 

% From gitik@math.tau.ac.il Thu Jul  6 01:55:37 1995
% From: Gitik Moti <gitik@math.tau.ac.il>
% Date: Thu, 6 Jul 1995 08:53:35 +0300
% To: mitchell@math.ufl.edu
% 
% 
% %This file is called def1200.tex you will need it
%to process
%the other file(s) you receive from me.

%-------------------------------------------------------
\def\today{\ifcase\month\or January\or February\or
March\or April\or May\or June\or July\or August\or
September\or October\or November\or December\fi
\space\number\day, \number\year}

%\font\lcmss=lcmssb8

%\font\cmss=cmssq10 at 10pt

%\font \rrrm=cmbx10 at 12pt

\def\dspace{\lineskip=2pt\baselineskip=18pt
\lineskiplimit=0pt}
\def\sspace{\lineskip=2pt\baselineskip=12pt
\lineskiplimit=0pt}

\font \bbrm=cmbx10 at 12pt

\font \ninerm= cmr10 at 9pt

\def\smalltype{\ninerm}

\def\bigtype{\bbrm}

\hsize=13.5cm
\magnification=1200
\def\ce{\centerline}

\def\hb{\hfill\break}

\def\title #1{\null\bigskip\ce{\bigtype #1}
\bigskip}
%for discretionary break in equation use

%Greek letters
\def\alp{\alpha}		
\def\bet{\beta}		
\def\gam{\gamma}		
\def\del{\delta}

\def\kap{\kappa}
\def\lam{\lambda}		
\def\sig{\sigma}		

\def\ome{\omega}		

%Caligraphic roman letters

\def\calC{{\cal C}}

\def\calK{{\cal K}}

\def\calP{{\cal P}}

%Bold roman letters

%bold Greek capital letters

%Capital roman double letters
%\def\CC{{\rlap {\raise 0.4ex \hbox{$\scriptscriptstyle
%|$}}
%\hskip -0.15em C}}
%\def\AA{{\mathchoice
%{I\hskip -3.7pt {\rm A}}
%{I\hskip -3.7pt {\rm A}}
%{I\hskip -3.1pt {\rm A}}
%{I\hskip -2.5pt {\rm A}}}}
%\def\BB{{I\!\!B}}
%\def\EE{{I\!\!E}}
%\def\FF{{I\!\!F}}
%\def\NN{{I\!\!N}}
%\def\PP{{I\hskip-2.5pt P}}
%\def\QQ{{\rlap {\raise 0.4ex
%\hbox{$\scriptscriptlsubspacesstyle |$}}
%\hskip -0.1em Q}}
%\def\RR{{I\!\!R}}
%\def\ZZ{{Z\!\!\! Z}}
    
%Special fonts
\font\tenboldgreek=cmmib10
 \font\sevenboldgreek=cmmib10 at 7pt
\font\fiveboldgreek=cmmib10 at 7pt
\newfam\bgfam
\textfont\bgfam=\tenboldgreek
\scriptfont\bgfam=\sevenboldgreek
\scriptscriptfont\bgfam=\fiveboldgreek

%\def\bg{\fam6}
%\mathchardef\alpha="700B
%\def\bfalp{{\fam=\bgfam\balp}}
\mathchardef\ggarrow="7010

\font\tengerman=eufm10 \font\sevengerman=eufm7
\font\fivegerman=eufm5
\font\tendouble=msym10 \font\sevendouble=msym7
\font\fivedouble=msym5

\textfont4=\tengerman \scriptfont4=\sevengerman
\scriptscriptfont4=\fivegerman
\newfam\dbfam
\textfont\dbfam=\tendouble \scriptfont\dbfam=
\sevendouble
\scriptscriptfont\dbfam=\fivedouble

\mathchardef\ng="702D
\mathchardef\dbA="7041
\mathchardef\sm="7072
\mathchardef\nvdash="7030
\mathchardef\nldash="7031
\mathchardef\lne="7008
\mathchardef\sneq="7024
\mathchardef\spneq="7025
\mathchardef\sne="7028
\mathchardef\spne="7029
\mathchardef\ltms="706E
\mathchardef\tmsl="706F

\mathchardef\dbA="7041

%Euler Fraktur letters

%Capital roman double letters
\mathchardef\dbA="7041 
\mathchardef\dbB="7042 
\mathchardef\dbC="7043 
\mathchardef\dbD="7044 
\mathchardef\dbE="7045 
\mathchardef\dbF="7046 
\mathchardef\dbG="7047 
\mathchardef\dbH="7048 
\mathchardef\dbI="7049 
\mathchardef\dbJ="704A 
\mathchardef\dbK="704B 
\mathchardef\dbL="704C 
\mathchardef\dbM="704D 
\mathchardef\dbN="704E 
\mathchardef\dbO="704F 
\mathchardef\dbP="7050 
\mathchardef\dbQ="7051 
\mathchardef\dbR="7052 
\mathchardef\dbS="7053 
\mathchardef\dbT="7054 
\mathchardef\dbU="7055 
\mathchardef\dbV="7056 
\mathchardef\dbW="7057 
\mathchardef\dbX="7058 
\mathchardef\dbY="7059 
\mathchardef\dbZ="705A 

\def\nek{,\ldots,}
\def\sdp{\times \hskip -0.3em {\raise 0.3ex
\hbox{$\scriptscriptstyle |$}}} % semidirect product

%words in roman font

\def\cf{{\rm \,cf\,}}

\def\dom{\mathop{\rm dom}\nolimits}

%\def\rank{\rm rank}

%overlined math alphabet

\def\op{{\overline p}}

%overlined Greek alphabet

%underlined math alphabet

%underline Greek alphabet

%math alphabet with hat

%Greek alphabet with hat

%roman with widetilde

%Greek alphabet with widetilde

\def\ddownarrow{\big\downarrow \hskip-0.70em\raise
2pt\hbox {$\big\downarrow$}}
\def\longright #1#2 {\smash{\mathop{\hbox to
#1pt {\rightarrowfill}}\limits_{#2}}}
\def\sqr#1#2{{\vcenter{\hrule height.#2pt\hbox{\vrule
width.#2pt height#1pt \kern#1pt \vrule width.#2pt}
\hrule height.#2pt}}}

\def\buildrul#1\under#2{\mathrel{\mathop{\null#2}
\limits_{#1}}}

\def\boxit#1{\vbox{\hrule\hbox{\vrule\kern3pt
\vbox{\kern3pt#1 \kern3pt}\kern3pt\vrule}\hrule}}

\def\prodl{\prod\limits}

\def\subheading#1{\medskip\goodbreak\noindent{\bf
#1.}\quad}

\def\sect#1{\goodbreak\bigskip\centerline{\bf#1}
\medskip}
\def\pr{\smallskip\noindent{\bf Proof:\quad}}
\def\onumber #1{\ooalign{\hfil\raise.07ex\hbox{
\hfill$\scriptstyle \,#1$\hfil}
\cr\cr{$\bigcirc$}}}
\def\onumber c{\ooalign{\hfil\raise.07ex\hbox
{\hfill$\scriptstyle \,c$\hfil}
\cr\cr{$\bigcirc$}}}
\def\alpcirc {\ooalign{\hfil\raise.07ex
\hbox{\hfill$\scriptstyle\alp\;$\hfill}\cr\cr
{$\bigcirc$}}}
\def\astcirc {\ooalign{\hfil\raise.07ex
\hbox{\hfill$\textstyle\ast\;$\hfill}\cr\cr
{$\bigcirc$}}}

\def\longmapright #1#2 {\smash{\mathop{\hbox to
#1pt {\rightarrowfill}}\limits^{#2}}}
\def\longmapleft #1 #2 {\smash{\mathop{\hbox to
#1 pt {\leftarrowfill}}\limits^{#2}}}

\def\references#1{\goodbreak\bigskip\par\centerline
{\bf References}\medskip\parindent=#1pt}
\def\ref#1{\par\smallskip\hang\indent\llap{\hbox
to \parindent{#1\hfil\enspace}}\ignorespaces}

\def\back{{\raise 2.5pt\hbox{$\,\scriptscriptstyle
\backslash\,$}}}
\def\bks{{\backslash}}
\def\part{\partial}
\def\lwr #1{\lower 5pt\hbox{$#1$}\hskip -3pt}
\def\rse #1{\hskip -3pt\raise 5pt\hbox{$#1$}}
\def\lwrs #1{\lower 4pt\hbox{$\scriptstyle #1$}
\hskip -2pt}
\def\rses #1{\hskip -2pt\raise 3pt\hbox
{$\scriptstyle #1$}}

\def\<#1{\left\langle{#1}\right\rangle}

\def\subinbn{{\subset\hskip-8pt\raise 0.95pt
\hbox{$\scriptscriptstyle\subset$}}}

\def\llvdash{\mathop{\|\hskip-2pt
\raise 3pt\hbox{\vrule height 0.25pt width 1.5cm}}}

\def\lvdash{\mathop{|\hskip-2pt \raise 3pt\hbox
{\vrule height 0.25pt width 1.5cm}}}

\def\fakebold#1{\leavevmode\setbox0=\hbox{#1}%
  \kern-.025em\copy0 \kern-\wd0
  \kern .025em\copy0 \kern-\wd0
  \kern-.025em\raise.0333em\box0 }

\font\msxmten=msxm10
\font\msxmseven=msxm7
\font\msxmfive=msxm5
\newfam\myfam
\textfont\myfam=\msxmten
\scriptfont\myfam=\msxmseven
\scriptscriptfont\myfam=\msxmfive
\mathchardef\rhookupone="7016
\mathchardef\ldh="700D
\mathchardef\leg="7053
\mathchardef\ANG="705E
\mathchardef\lcu="7070
\mathchardef\rcu="7071
\mathchardef\leseq="7035
\mathchardef\qeeg="703D
\mathchardef\qeel="7036
\mathchardef\blackbox="7004
\mathchardef\bbx="7003
\mathchardef\simsucc="7025

\def\rhookup{{\fam=\myfam \rhookupone}}

\def\bigsquare{{\fam=\myfam\bbx}}

\font\tencaps=cmcsc10
\def\smallcaps{\tencaps}

\def\author#1{\bigskip\ce{\smallcaps #1}\medskip}

\def\tagg{^{\prime\prime}}

\def\upddots{\mathinner{\mkern
1mu\raise 1pt \hbox{.}\mkern 2mu \mkern
2mu \raise 4pt\hbox{.}\mkern 1mu \raise 7pt\vbox
{\kern 7 pt\hbox{.}}} }

\def\varchi{\ooalign{{\raise
1.385pt\hbox{$\chi$}}\crcr\hbox{--}\crcr}}

\def\trianarrow{{\raise 2pt\hbox to 0.50cm
{\hrulefill}\triangleright}}
\def\Chi{{\raise 3pt\hbox{$\chi$}}}

\font\b=cmr10 scaled \magstep4

\def\bigzerou{\smash{\lower1.7ex\hbox{\b 0}}}
\def\bigast{\smash{\lower1.7ex\hbox{\b *}}}

\def\leaderfill{leaders\hbox to 5em{\hss\hss}\hfill}
\newcount\notenumber
       
       \def\note#1{\advance\notenumber by
1\footnote{$^{\the\notenumber}$} {\sspace\smalltype #1}}

\null
\overfullrule=0pt
{\nopagenumbers
\sect{ON CLOSED UNBOUNDED SETS CONSISTING}
\vskip -.5truecm
\sect{OF FORMER REGULARS}
\dspace
\sect{Moti Gitik}
\bigskip
\ce{School of Mathematical Sciences}
\ce{Sackler Faculty of Exact Sciences}
\ce{Tel Aviv University}
\ce{Ramat-Aviv, 69978 Israel}
\bigskip
{\narrower\medskip{\narrower\medskip
\sect{Abstract}

A method of iteration of Prikry type forcing
notions as well as a forcing for adding clubs is
presented.  It is applied to construct a model
with a measurable cardinal containing a club of
former regulars, starting with $o(\kap)=\kap +1$.  On
the other hand, it is shown that the
strength of above is at least $o(\kap)=\kap$.\medskip}
\medskip}

\vfill\eject}
\count0=1
\null
\dspace
\def\llvdash{\mathop{\|\hskip-2pt \raise 2pt\hbox{\vrule
height 0.25pt width 0.50cm}}}
\def\sit{{\buildrul\sim\under t}} 
\def\sima{{\buildrul\sim\under a}} 
\def\simB{{\buildrul\sim\under B}} 
\def\simD{{\buildrul\sim\under D}} 
\def\simP{{\buildrul\sim\under P}} 
\def\simX{{\buildrul\sim\under X}} 
\def\simcalP{{\buildrul\sim\under\calP}} 
\def\simapr{{\vtop{\offinterlineskip\hbox{$a'$}
\hbox{$\sim$}}}}
\def\simda{{\vtop{\offinterlineskip\hbox{$D_\alp$}
\hbox{$\sim$}}}}
\def\simai{{\vtop{\offinterlineskip\hbox{$a_i$}
\hbox{$\sim$}}}}
\def\simao{{\vtop{\offinterlineskip\hbox{$a_1$}
\hbox{$\sim$}}}}
\def\simat{{\vtop{\offinterlineskip\hbox{$a_2$}
\hbox{$\sim$}}}}
\def\simaz{{\vtop{\offinterlineskip\hbox{$a_0$}
\hbox{$\sim$}}}}

\def\simpstar{{\vtop{\offinterlineskip\hbox{$p^*$}
\hbox{$\sim$}}}}
\def\simpb{{\vtop{\offinterlineskip\hbox{$p_\bet$}
\hbox{$\sim$}}}}
\def\simpbstar{{\vtop{\offinterlineskip\hbox{$p_\bet^*$}
\hbox{$\sim$}}}}
\def\simpbstarr{{\vtop{\offinterlineskip
\hbox{$p_\bet^{**}$} \hbox{$\sim$}}}}
\def\simpg{{\vtop{\offinterlineskip\hbox{$p_\gam$}
\hbox{$\sim$}}}}
\def\simppri{{\vtop{\offinterlineskip\hbox
{$p_{i_1\gam_1}'$}\hbox{$\sim$}}}}
\def\simpbeto{{\vtop{\offinterlineskip\hbox{$p_{\bet_1}$}
\hbox{$\sim$}}}}
\def\simpbetostar{{\vtop{\offinterlineskip
\hbox{$p^*_{\bet_1}$}\hbox{$\sim$}}}}

\def\simpid{{\vtop{\baselineskip=5pt\hbox{$p_{i\del}$}
\hbox{$\sim$}}}}
\def\simpin{{\vtop{\baselineskip=5pt\hbox{$p_{i\nu}$}
\hbox{$\sim$}}}}
\def\simpgpr{{\vtop{\offinterlineskip\hbox{$p_\gam'$}
\hbox{$\sim$}}}}
\def\simpgprn{{\vtop{\offinterlineskip\hbox{$p_{\gam_n}'$}
\hbox{$\sim$}}}}
\def\simpgstar{{\vtop{\offinterlineskip\hbox{$p^*_\gam$}
\hbox{$\sim$}}}}
\def\simpjd{{\vtop{\offinterlineskip\hbox{$p_{j\del}$}
\hbox{$\sim$}}}}
\def\simqg{{\vtop{\offinterlineskip\hbox{$q_\gam$}
\hbox{$\sim$}}}}
\def\simQa{{\vtop{\offinterlineskip\hbox{$Q_\alp$}
\hbox{$\sim$}}}}
\def\simQg{{\vtop{\offinterlineskip\hbox{$Q_\gam$}
\hbox{$\sim$}}}}
Suppose that $\kap$ is an inaccessible cardinal.  We wish
to find a generic extension (usually cardinal preserving)
such that $\{\alp < \kappa \mid \alp$ is regular
in $V\}$ contains a club.  Radin [Ra]
introduced a basic method to do
this.  Simply start with a measurable $\kappa$ with
$o(\kappa) = \kappa^+$ and then force with the Radin
forcing constructed from $o(\kappa) = \kappa^+$. 
If one wishes to keep $\kappa$ a measurable in the
extension, then a weak repeat point suffices.
Both facts are proved by Mitchell [Mi1],

We show how to reduce assumptions rendering the above
possible.

A method of iteration generalizing those of
[Gi1] is presented.  Then a variant of it is
used to iterate forcing for shooting clubs.  We
think that this method of iteration can be
applied to other distributive forcings as well.   

We like to thank the referee for pointing out
that the proof of section 2 gives only
$o(\kap)=\kap$ and not $o(\kap)\kap +1$  as was
claimed in the previous version, for long and
detailed list of corrections and for his
requests on explaining certain parts of the paper.

\sect{1.~~A Forcing Construction}

We will now prove the following
two theorems.
\proclaim Theorem 1.1.  Suppose that $\kappa$ is
an inaccessible cardinal such that for
every $\del < \kappa$ the set of $\alp$'s below
$\kappa$ with $o(\alp) \ge \del$ is stationary.
Then there is a cardinal preserving generic
extension such that the set $\{\alp <
\kappa\mid \alp$ is regular in $V\}$ contains a
club.

\proclaim Theorem 1.2.  Suppose that $\kappa$ is
a measurable cardinal with $o(\kappa) = \kappa +
1$.  Then there is a cardinal preserving
extension satisfying the following:
\itemitem{(1)} $\kappa$ is a measurable,
\itemitem{(2)} $\{\alp < \kappa \mid \alp$ is a
regular in $V$ contains a club$\}$.

The proofs of these theorems use an
iteration Prikry type forcing notion that
was introduced in [Gi1].

Basically for every $\alp < \kappa$ with $\alp >o
(\alp )>o$ we are forcing a Prikry or Magidor sequence to
$\alp$  without adding new bounded subsets.  The order
type of the sequence is $\ome^{o(\alp)}$,  where the
exponentiation is the ordinal one.  The forcing used for
this is $\langle\calP(\alp,o(\alp)),\le, \le^*\rangle$
which was introduced in [Gi1]. The definition of it is
quite long and will not be used explicitly here.  The
only facts used in connection with this are stated below.

\subheading{Fact 1} $\langle\calP(\alp,o(\alp))$,
$\le$, $\le^*\rangle$  satisfies Prikry
condition and it is $\alp$-weakly closed, i.e.
\item{(a)} $\le\ \supseteq\ \le^*$;
\item{(b)} for every $p\in\calP(\alp,o(\alp))$  and every
statement $\sig$  of the forcing $\langle\calP
(\alp,o(\alp)),\le \rangle$ there is $q\ge^{\!\!\!{}^*}p$
 deciding $\sig$; (c) $\le^*$  is $\alp$-closed.

\subheading{Fact 2} $\langle \calP
(\alp,o(\alp)),\le \rangle$  satisfies
$\alp^+$-c.c. and has cardinality $2^\alp$.  

\subheading{Fact 3} $\langle\calP(\alp,
o(\alp)),\le \rangle$  adds a closed cofinal in
$\alp$  sequence consisting of regular in $V$
cardinals of order type $\ome^{o(\alp)}$ and
almost containing in every club of $\alp$ of $V$.

We are going to iterate $\calP(\alp,0(\alp))$'s
using the iteration of [Gi1]. For the benefit of
the reader let us give a precise definition.

Let $A$  be a set consisting of $\alp$'s such
that $\alp <\kap$  and $\alp >o(\alp)>0$.
Denote by $A^\ell$  the closure of the set $\{\alp
+1\mid \alp\in A\}\cup A$.  For every $\alp\in
A^\ell$  define by induction $\calP_\alp$ to be
the set of all elements $p$  of the form
$\langle\simpg\mid\gam\in g \rangle$  where 
\item{(1)} $g$ is a subset of $\alp\cap A$.
\item{(2)} $g$ has an Easton support, i.e. for
every inaccessible $\bet\le\alp$,  $\bet
>|\dom g\cap\bet |$; 
\item{(3)}
$p\rhookup\gam =\langle\simpb\mid\bet\in\gam\cap g
\rangle\in\calP_\gam$  and
$p\rhookup\gam\llvdash_{\calP_\gam}$
``$\simpg\in\calP (\gam,o(\gam))$" for every
$\gam\in\dom g$

Let $p=\langle\simpg\mid\gam\in g\rangle\ ,\ q=\langle
\simqg\mid\gam\in f\rangle$ be elements
of $\calP_\alp$.  Then $p\ge q$  ($p$  is
stronger than $q$) if the following holds: 
\item{(1)} $g\supseteq f$
\item{(2)} for every $\gam\in f$ 
$p\rhookup \gam\llvdash_{\calP_\gam}$
``$\simpg\ge\simqg$  in the forcing
$\calP(\gam,0(\gam))$"
\item{(3)} there exists a finite subset $b$ of
$f$ so that for every $\gam\in f\bks b$,
$p\rhookup \gam\llvdash$ ``$\simpg\ge^*\simqg$
in $\calP(\gam,0(\gam))$".

If $b=\emptyset$,  then let $p\ge^*q$.

By [Gi1], such iteration preserves cardinals.
The final forcing $\calP_\kappa$ satisfies $\kappa$ -
c.c.  We refer to [Gi1] for more details.

\subheading{Proof of Theorem 1}
Let $G_\kappa$ be a generic subset of $\calP_\kappa$.
We force over $V [G_\kappa]$ with
the forcing $P[E] = \{ d|d$ is closed, bounded
subset of $E\}$ ordered by the end extension, where
$$E = \{ \alp < \kappa\mid \alp\ \hbox{is a regular in}
\ V \} \ ,$$
i.e. the usual forcing for adding a club via a stationary
set. 

By Avraham-Shelah [Av-Sh], if $E$  is fat then
this forcing is $(\kap ,\infty)$-distributive.
Where fatness means that for every $\del <\kap$
and every closed unbounded subset $C$  of $\kap$
there is a closed subset $s\subseteq C\cap E$ of
order type $\del$. 

So it is enough to show that $E$  is fat in
$V[G_\kap]$.  Let $\del <\kap$ and
$C\subseteq\kap$  be a club in $V[G_\kap]$.  By
$\kap$-c.c. of the forcing $\calP_\kap$ there is
$C^*\subseteq C$  a club which belongs to $V$.
Then pick an $\alp^*$  a limit point of
$C^*\cap\{\alp <\kap\mid o(\alp)\ge\del\}$. 
By Fact 3. there will be a
closed subset of $C^*\cap\alp^*\cap E$  of order
type $\del$. This completes the proof.\hfill$\bigsquare$ 

We now turn to the proof of Theorem 1.2 which
takes up the rest of this section.  

\subheading{Proof of Theorem 1.2}
Denote by $U$ the measure concentrating on
$\alp$'s with $o(\alp) = \alp$. 

We preserve the notation of Theorem 1.1. 

Let us explain the idea of the proof.  We like
to force with $P[E]$ a club through the set $E$
of regular cardinals below $\kap$.  So lots of
changes of cofinalities are needed.  An
additional task is to preserve measurability of
$\kap$.  So for a set of $\alp$'s in $U$  we
need to use $P[E\cap\alp]$.  The problem with it
is that this forcing is not closed or weakly
closed but only distributive.  Iteration of
$P[E\cap\alp]$'s in Easton fashion even $\ome$-times
collapses cardinals.  The key idea will be to
embed $P[E\cap\alp]$  into $\calP(\alp,\bet)$  $*$ 
$\bet$-closed, where $\bet <\alp$  and we will allow to
change $\bet$ from time to time in order to increase
the degree of closeness.  Then the scheme of iteration
of [Gi1] will be generalized to make possible
the iteration of less closed forcing notions.

%\proclaim Lemma 1.3. Let $\bet < 0(\alp)$. Then
%$P[E\cap\alp]$  embeds into
%$\calP(\alp,\bet)* \bet$-closed.

%\subheading{Remark}  Actually $P[E\cap\alp]$
%embeds trivially into $\calP(\alp,\bet)*P[E\cap\alp]$
%and $P[E\cap\alp]$  in $V^{\calP(\alp,\bet)}$ contains
%a dense $\bet$-closed set.  

%\pr Let $b$ be the closed cofinal sequence in
%$\alp$  introduced by $\calP(\alp,\bet)$.  The
%set of all $p\in P[E\cap\alp]$  such that $\sup (p)\in b$
%is dense in $P[E\cap\alp]$.  Since all elements of $b$ 
%are regular in $V$,  we are free to take unions of
%conditions inside this dense set.\hfill$\bigsquare$

Let us start with a generalization of the iteration
process of [Gi1].

We describe a general scheme of iteration.  A
specific forcing used for the proof of 1.2 will
be defined later.

Let $A$  be a set consisting of inaccessible
cardinals.  Denote by  $A^\ell$  the closure
of the set $A\cup \{\alp+1\mid\alp\in A\}$.  We
define an iteration $\langle\calP_\alp$,
$\simQa\mid \alp\in A^\ell\rangle$.
For every $\alp\in A^\ell$  define by induction
$\calP_\alp$  to be the set of all elements $p$
of the form $\langle\simpg\mid\gam\in g\rangle$, where
\item{(1)} $g$ is a subset of $\alp\cap A$; 
\item{(2)} $g$  has an Easton support, i.e. for
every inaccessible $\bet\le\alp\quad\bet >|\dom
g\cap\bet |$;
\item{(3)} for every $\gam\in\dom g$
$$p\rhookup\gam =\langle\ \simpb|\bet
<\gam\rangle\in\calP_\gam$$ 
and
$p\rhookup\gam\llvdash_{\calP_\gam}$
``$\simpg$  is a condition in
the forcing $\langle\simQg,
\le_\gam, \le^*_\gam\rangle$ satisfying the Prikry
condition and of cardinality below the next element
of $\check A$ above $\check\gam$".

The difference here is that we do not require
that $\langle\simQg, \le^*_\gam\rangle$ is $\gam$-closed.

Let $p=\langle\simpg\mid\gam\in
g\rangle$  and $q=\langle\simqg |\gam\in f\rangle$ be
elements of $\calP_\alp$.  Then $p\ge q$ ($p$ is stronger
than $q$) if the following holds:
\item{(1)} $g\supseteq f$
\item{(2)} for every $\gam\in f$

$p\rhookup \gam\llvdash_{\calP_\gam}$ ``$\simpg\ge_\gam
\simqg$  in the forcing $\simQg$"
\item{(3)} there exists a finite subset $b$  of
$f$ so that for every $\gam\in f\bks b$
$p\rhookup\gam\llvdash_{\calP_\gam}$
``$\simpg\ge_\gam^{\!{}*}\simqg$ in the forcing $\simQg$".

If the set $b$  in (3) is empty we call $p$ a
direct extension of $q$  and denote this by
$p\ge^{\!{}*}q$.

\proclaim Lemma 1.3.  Let $\alp\in A^\ell$,
$p\in \calP_\alp$  and $\sig$  be a statement of
the forcing language appropriate for
$\calP_\alp$.  Then there is a direct extension
$p^*$  of $p$  deciding $\sig$.

\pr Suppose otherwise.  Let $p=\langle\simpg\mid\gam\in
g\rangle$
and $\bet$  be the minimal element of $g$.  We
assume that $g\not=\emptyset$, otherwise any
extension of it is direct.

Let $G$  be a generic subset of $\calP_{\bet +1}$,  
so that $p\rhookup \bet +1\in G$.  We shall mean
by $p\bks (\bet +1)=\langle 
\simpg\mid \gam\in g\bks (\bet +1)\rangle$ the
interpretation of it in $V[G]$, i.e. an element of the
forcing $\calP_\alp/G$.  Define now $p^*\in\calP_\alp/G$.
If there exists some $q\in\calP_\alp/G$  a
direct extension of $p\bks (\bet +1)$  deciding
$\sig$,  then set $p^*$ to some such $q$.
Otherwise, set $p^*=p\bks (\bet +1)$. 
Let $p^*=\langle\simpgstar
|\gam\in g^*\rangle$.  Then $g^*\in V[G]$.  But, since
$\calP_{\bet +1}$  is of small cardinality (see
(3) of the definition of $\calP_\alp$) $g^*$
can be easily replaced by a set in $V$.  Suppose
that is already the case.  Then $p\rhookup (\bet
+1)\cup p^*$  will be a condition in
$\calP_\alp$  and it will be a direct extension
of $p$.

\subheading{Claim 1.4.1} There is 
$\simpbstar$  such
that $\phi \llvdash_{\calP_\bet}$ ``$\simpstar\ge
\simpbstar$" and
$\langle p^*_\bet\rangle\llvdash_{\calP_{\bet +1}}$
$$\hbox{``}p^*=p\bks (\bet +1)\hbox{"}$$

\pr Let $G\subseteq\calP_\bet$  be generic.
Work in $V[G]$.  Let $p^*_\bet$  be a direct
extension of $p_\bet$  (in $Q_\bet$)  deciding
the statement ``$p^*=p\bks (\bet +1)$".  Suppose
for a moment that the decision is negative. Then
$p^*_\bet \llvdash_{Q_\bet}$  ``$p^*$  is a
direct extension of $p\bks (\bet +1)$  deciding
$\sig$", by the choice of $p^*$.  Now pick a
direct extension $p^{**}_\bet$  of $p^*_\bet$  so
that   
$$p^{**}_\bet \llvdash_{Q_\bet}\hbox{``} p^*\llvdash
{}^{\!i}\sig\hbox{"}\ ,$$
where ${}^{\!\circ}\sig =\sig$  and  ${}^{\!1}
\sig={}^{\!\neg}\sig$.  Pick some $r\in G$  so
that 
$$r^\cap\langle\simpbstarr\rangle\llvdash_{\calP_{\bet
+1}}\hbox{``}
p^*\llvdash{}^{\!i}\sig\hbox{"}\ .$$
Then $r^\cap\langle\simpbstarr\rangle^\cap p^*$  will be
a direct extension of $p$ forcing ${}^{\!i}\sig$.  Which
contradicts our assumption.
\hfill$\bigsquare$ of the claim.\break

Now we replace $\simpb$ in $p$  by $\simpbstar$.
Denote the resulting condition by $p(\{ \bet \})$.
Let $\bet_1$  be the second element of $g$.  We
proceed as above replacing $\bet$  by $\bet_1$
and $p$  by $p(\{\bet \})$.  This will define 
$\simpbeto$.  Then
set $p(\{\bet_0,\bet_1\})=\langle\simpbstar,\simpbetostar
\rangle^\cap\langle\simpg\mid \gam\in g\bks\{\bet,\bet_1
\}\rangle$. Continue in the
same fashion.  Finally, after going through all
the elements of $g$  we will obtain a condition
$p(g)=\langle\simpgstar\mid\gam\in g\rangle$ which is a
direct extension of $p$.  Now let $q\ge p(g)$
be a condition deciding $\sig$.  By the
definition of extension, there is a maximal
$\gam\in g$  such that $q\rhookup
\gam\llvdash_{\calP_\gam}$ ``$\simqg$ is not a direct
extension of $\simpg$".  But this
will contradict the choice of $p(g\cap (\gam
+1))$, since above $\gam$, $q\bks\gam$  will be
a direct extension of $p\bks\gam$    deciding $\sig$. 
Contradiction.\hfill$\bigsquare$

Now let us define the iteration needed for the
proof of 1.2. 

Let $\alp <\kap$ be an ordinal with
$o(\alp)=\alp$.  We like to define $Q_\alp$.
But first, let us consider the forcing
$\calP(\alp,\bet)*P[E\cap\alp]$,  where $\bet
<\alp$.  Define a $*$-ordering on
$\calP(\alp,\bet)*P[E\cap\alp]$  by setting $p=\langle
b, \sima\rangle\le^*p'=\langle
b',\simapr\rangle$  if $p=p'$
or (i) $p\le p'$ and (ii) $b\le^*b'$  in
$\calP(\alp,\bet)$.

\proclaim Lemma 1.4.  The forcing $\langle
\calP(\alp,\bet)*P[E\cap\alp],\le,\le^*\rangle$
satisfies the Prikry condition. 

The proof follows easily from Fact 1.

\proclaim Lemma 1.5.0.  Every condition in
$\calP(\alp,\bet)*P[E\cap\alp]$  can be
extended to a condition of the form $\langle
b,\check c\rangle$. 

It follows by Fact 1 since $\calP(\alp,\bet)$  does not
add new bounded subsets to $\alp$.

Let $D_{\alp,\bet}=\{\langle b,\sima\rangle\in
\calP(\alp,\bet)*P[E\cap\alp]\mid
b\llvdash_{\calP(\alp,\bet)}\hbox{``}\sup\sima\in\simB$",
where $\simB$  is the canonical name of the closed
cofinal in $\alp$ sequence added by
$\calP(\alp,\bet)\}$.

\proclaim Lemma 1.6.  For every
$p\in\calP(\alp,\bet)*P[E\cap\alp]$  there exists
$q\in D_{\alp,\bet}$ such that $q\ge^*p$.

\pr Obviously, since $b\llvdash_{\calP(\alp,\bet)}\sup
\buildrul\sim\under a <\alp$,  where $p=\langle
b,\sima\rangle$.\hfill$\bigsquare$

\proclaim Lemma 1.7.  The order $\le^*$  is
$\bet$-closed over $D_{\alp,\bet}$  and hence
$\calP(\alp,\bet)*P[E\cap\alp]$  does not add
new bounded subsets of $\bet$.

\subheading{Remark} Notice that all cardinals
$\le\alp$ are collapsed to $\bet$,  since the
forcing with $P[E\cap\alp]$  over
$V^{\calP(\alp,\bet)}$  is isomorphic to the Levy
collapse ${\rm Col}(\bet,\alp)$.

\pr Let $\langle <b_i,\simai >\mid i<\tau <\bet\rangle$ 
be a
$\le^*$-increasing sequence of elements of
$D_{\alp,\bet}$.  Since
$\langle \calP(\alp,\bet),\le^*\rangle$  is
$\alp$-closed, there is $b\in\calP(\alp,\bet)$  such
that $b\ {}^*\!\!\ge b_i$  for all $i<\tau$.
Let $\sima$ be a name of the union of $a_i$'s.
Since for every $i<\tau$ $b\llvdash_{\calP(\alp,\bet)}
\tagg\sup(\simai)\in\buildrul\sim\under B\tagg$
where $\buildrul\sim\under B$  is a canonical
name of the generic closed cofinal in $\alp$
sequence, $b\llvdash_{\calP(\alp,\bet)}
\hbox{``}\sup\sima\in\simB$"
and in particular  $b\llvdash_{\calP(\alp,\bet)}\hbox{``}
\sup\sima$  is regular in $V$".  Hence
$\langle b,\sima \rangle$  is a condition, it
belongs to $D_{\alp,\bet}$  and it is $\le^*$ stronger
than each $\langle b_i,\simai\rangle\ (i<\tau)$.\hfill
$\bigsquare$

\subheading{Definition 1.8}  
\item{(1)} $Q_\alp=P[E\cap\alp]\cup
\bigcup_{\bet <\alp}(\calP(\alp,\bet)*P[E\cap\alp])$.
\item{(2)} The ordering of $Q_\alp$:
\item{(2a)} The ordering of $P[E\cap\alp]$  and
of $\calP(\alp,\bet)*P[E\cap\alp]$ $(\bet
<\alp)$  is the usual one
\item{(2b)} Let $c\in P[E\cap\alp]$  and $\langle
b,\sima\rangle\in\calP(\alp,\bet)*P[E\cap\alp]$,
for some $\bet <\alp$.
Define
$$\eqalign{\langle &b,\sima\rangle\ge c\ {\rm
iff}\ b\llvdash\sima\ge_{P[E\cap\alp]}\check c\cr
&c\ge \langle b,\sima\rangle\ {\rm iff}\
b\llvdash\check c\ge_{P[E\cap\alp]}\sima\cr}$$
\item{(2c)} Let $\langle b_i,\simai\rangle\in\calP
(\alp,\bet_i)*P[E\cap\alp]$
where $i=0,1,\bet_0\not= \bet_1<\alp$.
Define $\langle b_0,\simaz\rangle\le \langle
b_1,\simao\rangle$  iff there is $c\in
P[E\cap\alp]$  such that  $\langle b_0,\simaz\rangle\le
c\le\langle b_1,\simao\rangle$.

We also define $*$-ordering of $Q_\alp$. 

\subheading{Definition 1.9}  Let $p,q\in
Q_\alp$.  Set $p\le^* q$  iff $p=q$  or for some $\bet
<\alp$ $p=\langle b_1,\simao
\rangle$, $q=\langle b_2,\simat\rangle\in\calP(\alp,\bet)
*P[E\cap\alp]$  and (i) $p\le q$  (ii) $b_1\le^* b_2$  in
$\calP(\alp,\bet)$.

\proclaim Lemma 10.  (1) $P[E\cap\alp ]$  is
dense in $Q_\alp$.

(2) The forcing $\langle
Q_\alp,\le\rangle$  is equivalent to $\langle
P[E\cap\alp],\le\rangle$  and hence preserves
cofinalities of cardinals.

\pr $P[E\cap\alp]$  is dense in $Q_\alp$ by the
definition of the order.  Also 
conditions incompatible in $P[E\cap \alp]$ remaining
so in $Q_\alp$.\hfill$\bigsquare$ 

\proclaim Lemma 11.  Let $\bet <\alp$.  Suppose
that $\langle b,\sima\rangle\in
\calP(\alp,\bet)*P[E\cap\alp]$ and $\sig$ is a statement
of the forcing language of $Q_\alp$.  Then there is a
direct extension (in $Q_\alp$) of $\langle b,
\sima\rangle$  deciding $\sig$ in the forcing with
$Q_\alp$. 

\subheading{Remark}  The lemma actually shows the
Prikry condition above $\calP (\alp,\bet)*P[E\cap\alp]$ 
part of $Q_\alp$.  

The lemma will follow from the following
statement.

\proclaim Lemma 12.  For every $\bet <\alp$  $Q_\alp$ 
is a projection of $\calP(\alp,\bet)*P[E\cap \alp]$.

\pr It is enough to project a dense subset of
$\calP(\alp,\bet)*P[E\cap\alp]$ onto a dense
subset of $Q_\alp$.  Consider
$D\subseteq\calP(\alp,\bet)\times P[E\cap\alp ]$,  which
is dense in $\calP(\alp,\bet)*P[E\cap\alp]$.
There is such $D$ since $\calP(\alp,\bet)$  does not
add new bounded subsets of $\alp$.  For $\langle
b,a\rangle \in D$  set $\pi (\langle
b,a\rangle)=a$.  By Lemma 10, $rng\pi
=P[E\cap\alp]$ is dense in $Q_\alp$. It is trivial that
$\pi$ is a projection map.\hfill$\bigsquare$   

We are ready now to define the iteration.  The
definition will be as above only $\bet$'s for
$Q_\alp$ will be picked generically more carefully. 
This is needed for cardinals preservation.

Let $A$  and $A^\ell$  be as above.  For every
$\alp\in A^\ell$  we define $\calP_\alp$  by
induction.

\subheading{Definition 13}  A forcing notion
$\calP_\alp$  consists of all elements $p$ of
the form $\langle\buildrul\sim\under p\!\!{}_\gam\mid
\gam\in g\rangle$  where
\item{(1)} $g$  is a subset of $\alp\cap A$;
\item{(2)} $g$  has an Easton support;
\item{(3)} for every $\gam\in \dom g$
$p\rhookup\gam =\langle\simpb\mid\bet <\gam\rangle\in
\calP_\gam$ and
$p\rhookup\gam\llvdash_{\calP_\gam}\tagg\simpg$  is a
condition in
either the forcing $\simQg$,
if $o(\gam)=\gam$ or in the forcing
$\calP(\gam,o(\gam))$,  if $o(\gam)<\gam$"
\item{(4)} for every $\tau\le\alp$ the set
$\{\gam<\alp\mid\simpg\in
\buildrul\sim\under P[E\cap\gam]\}$ or $\gam\ge \tau$ 
and for some $\bet <\tau$  $\simpg\in\buildrul\sim\under
\calP(\gam,\bet)*\buildrul\sim\under P[E\cap
\gam])\}$  is finite.

The ordering on $\calP_\alp$  is defined without
changes.

The definition and Lemma 1.10 insures that for
every $\gam$  with $o(\gam)=\gam$  the actual
forcing used over $\gam$ is $P[E\cap\gam]$.  But
in every separate condition $p\in\calP_\alp$
only finitely many $\gam$'s with $p_\gam\in
P[E\cap\gam]$  are allowed (the condition (4)).
The reason for this is to insure that the
iteration preserves cardinals.  Intuitively,
finite iteration of forcings $P[E\cap\gam]$'s
does no harm.  In order to do infinite iterations
(even of the length $\ome$), we like to have in
advance some information about closed pieces of
$E$.  Forcings $\calP (\gam,\bet)$'s are
actually used for this purpose.  Namely canonical
generic sequences produced by such forcings.

\proclaim Lemma 14.  $\langle
\calP_\alp,\le,\le^*\rangle$  satisfies the
Prikry condition.

The proof repeats the proof of Lemma 13.  The
additional condition (4) has no effect on it.

\proclaim Lemma 15.  Let $\gam <\alp$.
Suppose that $G$  is a generic subset of
$\calP_{\gam +1}$.  Then the forcing $\calP_\alp/G$
does not add new subsets of $\gam$. 

\pr Suppose that some $p\llvdash\buildrul\sim\under a
\subseteq\check\gam$.  It is enough to show the following.

\subheading{Claim 1.15.1}  There are $b\in V[G]$
and $\op\ge p$ such that
$\op\llvdash\buildrul\sim\under a=\check b$.

\subheading{Remark} $\op$ is not required to
be a direct extension of $p$. The reason for
this is the finite set of $\gam$'s in $p$
satisfying (4).

\pr Using Lemma 1.14, we define by induction a
$*$-increasing sequence $\langle p_i\mid
i\le\gam \rangle$  of extensions of $p$  so that
for each $i$ 
\item{(a)} $p_i$ decides the statement
``$i\in\sima$".
\item{(b)} if for some $\del$  with
$o(\del)=\del$  $\del\in g_i$,
then
$p_i\rhookup\del\llvdash_{\calP_\del}$
``$\simpid\in D_{\del,\bet}$   for some
$\bet >\gam$" where
$p_i=\langle\simpin\mid
\nu\in g_i\rangle$,
\item{(c)} if $j>i,\ j\le\gam$,  then for every
$\del\in g_j\bks g_i$  with $o(\del)=\del$
$$p_j\rhookup\del\llvdash\hbox{``}\simpjd
\in \simcalP(\del,\bet)*
\buildrul\sim\under P[E\cap\del]$$  
for some $\bet$  which is above sup$(g_i\cap\del)$".

First we extend $p$  to a condition $p'$
satisfying (b).  By (4) of Definition 1.13  it is
always possible.  But $p'$ need not be a direct
extension of $p$.  Now, by Lemma 1.9, find $p\tagg
\ge^*p'$  deciding ``$\check 0\in\sima$".  Let $p_0$ 
be an extension of $p\tagg$ obtained
as follows.  We replace each $\simpg$" in
$p\tagg$ for $\del$ with $0(\del)=\del$  which is
not in $p'$ by a stronger condition in
$\simcalP(\del,\bet)*\simP[E\cap \del]$ where $\bet$ is
picked to be
above every coordinate of $p'$  below $\del$.
By (3), (4) of Definition 1.13, only finitely
many coordinates $\del$  in $p\tagg$ should be
fixed this way.  So $p_0$ will be a condition
stronger than $p\tagg$ and a direct extension of
$p'$  deciding ``$\check 0\in
\sima$".

We continue by induction.  On successor stages we
proceed as above.  Suppose now that $i\le\gam$
is a limit ordinal and a sequence $\langle
p_\rho \mid\rho <i\rangle$  is defined and
satisfies the conditions (a) - (c) above.  Let
us argue that there is $p'$  a direct extension
of all of $p_\rho$'s $(\rho <i)$.  Let $p'$  be
obtained by taking direct extensions in each
coordinate separately.  This is possible by (b).
It is enough to show that such an obtained $p'$ is
a condition.  The only problematic point is (4)
of Definition 1.13.  By (b), $\{\gam <\alp\mid \simpgpr
\in\simP[E\cap\gam ]\}$  is empty.  So it remains to
show that for every $\tau <\alp$  the set
$\{\gam <\alp\mid\tau\le\gam$, for some $\bet
<\tau$  $\simpgpr\in\simcalP(\gam,\bet)*\simP[E\cap\gam
]\}$ is finite.  Suppose otherwise.
Let $\tau$, $\tau\le\gam_0<\gam_1<\cdots
<\gam_n<\cdots <\alp$ be witnessing this. Then
for each $n<\ome$  there is $\bet_n<\tau$  such
that $\simpgprn\in\simcalP(\gam_n,\bet_n)*\simP[E\cap
\gam_n]$.  For each $n<\ome$  let $i_n<i$  be
the least such that the coordinate $\gam_n$
appears in $p_{i_n}$.  Shrinking the set of 
indexes if necessary, we assume that the
sequence $\langle i_n\mid n<\ome \rangle$   
is strictly increasing.  But this is impossible
by (c). Since $\gam_1>\gam_0$, $\gam_1$  is in
$p_{i_1}$ and not in $p_{i_0}$,  $\gam_0$  is in
$p_{i_0}$ but $\simppri\in\calP (\gam_1,\bet_1)*\simP
[E\cap \gam_1]$  where $\bet_1 <\tau <\gam_0$. 
Contradiction.
So $p'$  is a condition.  Now we continue as in
successor stages.\hfill$\bigsquare$ of the claim.\break
\hfill$\bigsquare$

\proclaim Lemma 1.6.  $\calP_\kap$  satisfies
$\kap$-c.c. and preserves the cardinals.  

Follows from  Definition 1.13 and Lemma 1.15.

Let $G_\kap$ be a generic subset of
$\calP_\kap$.  Force with $P[E]$  over $V[G_\kap]$.
Let $C$ be a generic club.

\proclaim Lemma 1.17.  $\kap$  is a measurable
cardinal in $V[G_\kap,C]$.

\pr Let $U$ be the measure in $V$  concentrating
over $\{\alp <\kap\mid o(\alp)=\alp\}$.  Denote
by $i_u:V\to N_u\simeq Ult (V,U)$  the
corresponding elementary embedding. Since $\calP_\kap$
satisfies $\kap$-c.c. and using Claim 1.15.1 it is    
routine to extend $i_u$  (in $V[G_\kap, C])$ to
an embedding
$$i:V[G_\kap ]\longrightarrow N[G_{i_u(\kap)}]$$
where  $N$ is an iterated ultrapower of $N_u$
moving only ordinals in the interval
$(\kap^+,i_u(\kap))$.  We refer to [Gi1,3] or
[Gi-Sh] for such arguments.  

For the benefit of the reader, let us provide
more details on what is going on.

Working in $V[G_\kap,C]$  we define a normal
$V[G_\kap]$-ultrafilter over $\kap$  (i.e. a
normal ultrafilter over $\calP^{V[G_\kap]}(\kap)$).
Proceed as follows:  pick in $V$  an
enumeration $\langle {\vtop{\offinterlineskip\hbox
{$D_\alp$} \hbox{$\sim$}}}\mid\alp <\kap^+\rangle$ such
that in $N_U$\hb
$\emptyset\llvdash_{\calP_{i_u(\kap)}\big/
G_{\kap^*}C}\hbox{``}\simda$ is a $*$-dense open
subset of $i_u(\calP_\kap)/\calP_{\kap +1}$"
and if for some $\simD$  $\emptyset\llvdash_{\calP_{i_u
(\kap)}/G_\kap *C}$ ``$\simD$  is a $*$-dense
open subset of $i_u(\calP_\kap)/\calP_{\kap +1}$
then for some $\alp <\kap^+$ $\emptyset\llvdash\simD
\subseteq\simda$.

Then define a master condition sequence $\langle
r_\alp\mid \alp <\kap^+\rangle\in V[G_\kap,C]$ such
that 
\item{(a)} $r_\alp\le^*r_\bet$  for every $\alp\le\bet
<\kap^+$
\item{(b)} $r_\alp\in D_\alp [G_\kap, C]$  for
every $\alp <\kap^+$.
\item{(c)} $\langle r_\alp\mid \alp
<\bet\rangle\in N_u[G_\kap,C]$ for every $\bet
<\kap^+$. 
\item{}
\smallskip

Using $\langle r_\alp \mid \alp <\kap^+\rangle$
we define now $U^*\subseteq U$ by setting $X\in
U^*$ iff for some $p\in G_\kap *C$,  $\alp
<\kap^+$,  a $\calP_\kap$-name
$\simX$ of $X$, in $N_u$   
$$\langle p,r_\alp\rangle\llvdash_{i_u(\calP_\kap)}
\hbox{``}\check\kap\in\simX\hbox{"}\ .$$
It is easy to check that $U^*$  is a normal
$V[G_\kap]$-ultrafilter.  Moreover, since $P[E]$
does not add new sequences of the length less
than $\kap$  to $V[G_\kap]$,  $U^*$  is
$\kap$-complete in $V[G_\kap,C]$.  Hence
$Ult(V[G_\kap],U^*)$  is well-founded.  Let
$N^*$ be its transitive collapse and
$i:V[G_\kap]\longrightarrow N^*$ be the
corresponding elementary embedding.  By
elementarity $N^*$  is of the form
$N[G_{i(\kap)}]$  where $G_{i(\kap)}$  is
$N$-generic subset of $\calP_{i(\kap)}$.
Examining the structure of $N$  it is possible
to show that it is in fact an iterated
ultrapower of $N_u$  moving only ordinals in the
interval $(\kap^+,i_u(\kap))$.  But this is not
needed for further argument. 

Since $P[E]$  does not add new $<\kap$-sequences,
$N[G_{i(\kap)}]$ is closed under $<\kap$-sequences of
its elements. But really more is true:

\subheading{Claim 1.17.1} 
${}^\kap N[G_{i(\kap)}]\subseteq
N[G_{i(\kap)}]$.

\pr Let $\langle t_\alp\mid\alp <\kap\rangle$
be a sequence of elements of $N[G_{i(\kap)}]$.
Without loss of generality, we may assume that
each $t_\alp$  is an ordinal.  Let
$\buildrul\sim\under t$  be a canonical
$\calP_\kap *P[E]$  name of this sequence.  Then
for each $\alp$  $\buildrul\sim\under t(\alp)$
is a set of cardinality at most $\kap$  consisting
of pairs $\langle p,\check\del\rangle$  where
$p\in\calP_\kap *P[E]$  and $\del$  is an
ordinal.  For every $\alp,\gam <\kap$  such that
$o(\gam)=\gam$ let $\buildrul\sim\under t(\alp)
\rhookup\gam$  denotes the set of $\gam$  first
pairs $\langle p,\check f_\del(\gam)\rangle$ in
$\buildrul\sim\under t(\alp)$ such that
$p\in\calP_\gam *P[E\cap\gam ]$ where $f_\del$ 
represents $\del$  in $N_u$.  Set $\sit\rhookup\gam
=\langle\sit (\alp)\rhookup\gam\mid\alp
<\gam\rangle$.  Then, $i_u(\langle\sit\rhookup\gam\mid
\gam<\kap\rangle)(\kap)=\sit$.  Now define a
function $g\in V[G_\kap]$  representing $\langle
t_\alp\mid\alp <\kap\rangle$  in $N[i_u(\calP_\kap)]$.
Set $g(\gam)=\sit \rhookup\gam [G_\gam
*C_\gam]$,  where $G_\gam *C_\gam =G_{\gam +1}=
G_\kap\cap\calP_\gam *P[E\cap\gam]$.\hfill$\bigsquare$  

Now, in $V[G_\kap]$
there is a Rudin-Kiesler increasing commutative
sequence of ultrafilters $\langle U_\alp \mid
\alp <\kap \rangle$  over $\kap$.  Thus for each
$\alp <\kap$  the measure concentrating on
$\{\bet <\kap\mid o(\bet)=\alp\}$  in $V$
extends to $U_\alp$  in $V[G_\kap]$.  Actually,
such extensions are used to define
$\calP(\gam,\del)$'s.  Then, in $N[G_{i(\kap)}]$
we will have such a sequence over $i(\kap)$
of the length $i(\kap)$.  Form the direct
limit.  Let $k:N[G_{i(\kap)}]\to M[G_\lam]$
be the corresponding embedding, where $\lam=k
(i(\kap))$.  Let $j=k\circ i: V[G_\kap]\to
M[G_\lam]$.  Notice that since the length of the
sequence used to form the direct limit is $i(\kap)$  
and $\cf (i(\kap))=\kap^+$,  ${}^\kap
M[G_\lam]\subseteq M[G_\lam]$.  An additional
point here is that $\lam$  is a limit of
critical points of embeddings used in the direct
limit.  They are are singular cardinals in $M[G_\lam]$
but are regular in $M$.  Since $M$ is just an
iterated ultrapower of $\calK$,  by [Mi2].
Hence, in $V[G_\kap,C]$, $j(E)$  contains a club.  Now we
use this in the standard fashion to diagonalize
over $\kap^+$  dense subset of $(P[j(E)])^{M[G_\lam]}$.
It will produce a club $C^*\supseteq C$  which
is $M[G_\lam]$-generic. Notice that by Lemma
1.10(2), $Q_\kap$  is equivalent to $P[E\cap \kap]$.  So
$C$  is $Q_\kap$  generic.  So we obtain $j\subseteq
j^*:V[G_\kap,C]\to M[G_\lam,C^*]$. Hence $\kap$
is measurable in $V[G_\kap,C]$.\hfill$\bigsquare$   

\sect{2.~~On the strength of the existence of a club
of ``regulars"}

In this section we will show that the hypothesis used in
Theorems 1.1 are the best possible and
those of 1.2 are close to this.

The next theorem is basically due to Mitchell. 

\proclaim  Theorem 2.1.  Suppose that $\kappa$
is an inaccessible and the set $E = \{ \alp
< \kappa\mid \alp$ is regular in $\calK\}$ contains
a club.  Then for every $\del < \kappa$ the set
$\{\alp < \kappa\mid o^\calK (\alp) \ge \del\}$ is
stationary.

\pr  Let $\del < \kappa$. We show that $\{\alp <
\kappa\mid o^\calK (\alp) \ge \del\} = A_\del$ is
stationary.  Let $C$ be a club contained in $E$.
Choose some $\alp \in C$ of cofinality
$\del$.  Then $\cf^\calK\alp = \alp$ since $\alp
\in E$.  So its cofinality changed to $\del$.  Then by
Mitchell [Mi2], $o^\calK(\alp) \ge \del$.
Hence $C\cap A_\del\not=\emptyset$.\hfill$\bigsquare$

\proclaim Theorem 2.2.  Let $\kappa$ be a
measurable cardinal and a set $E =\{\alp <
\kappa\mid \alp\ \hbox{is regular
in}\ \calK\}$ contain a club.  Then
$o^\calK(\kappa)\ge \kap$.

\pr First note that by Mitchell [Mi3],
$\kappa^+ = (\kappa^+)^\calK$.  Let $U$ be
a normal measure over $\kappa$.  Consider
its elementary embedding $j_U: V \to N
\simeq Ult(V, U)$.  By Mitchell [Mi3],
$j_U\rhookup\enskip \calK = i$ is an iterated
ultrapower of $\calK$ by its measures.
Suppose that $o^\calK(\kappa)<\kappa$. Then 
$i(\kappa)$ is not a limit point of
iteration.  Since ${}^\kap N\subseteq N,\cf
i(\kap)\ge \kap^+$. 
That is, there is a last
measure in which the ultrapower reached
$i(\kappa)$, (or more precisely, the image
of its critical point).  Denote this
critical point by $\lam$.  Then for every
$\alp < i(\kappa)$ there is $f:[\kappa]^n
\to \kappa$ ($n < \omega )$ in $\calK$ and
$\kappa \le \beta_1 \nek \beta_n \le \lam$
such that\enskip $i(f) (\beta_1 \nek \beta_n) =
\alp$.

Using coding of $n$-tuples,  we can replace
$n$-placed functions by a 1-placed ones.
Then, for every $\alp < i(\kappa)$ (or
$g:\kappa \to \kappa)$ there will be $\beta
\le \lam$ and $f:\kappa \to \kappa$ in
$\calK$ such that $i(f) (\beta) = \alp$ (or
$i(f) (\beta) = j_U(g) (\kappa)$).

\proclaim Claim 1.  For every $\alp <
i(\kappa)$ there is $f \in {^\kappa \kappa}
\cap \calK$ such that $\alp \le i(f)
(\lam)$.

\pr Let $\alp < i(\kap)$ and $g \in {^\kappa
([\kap]^n)} \cap \calK$ be such that $i(g)
(\beta_1 \nek \beta_{n-1}, \lam) = \alp$,
where $\beta_1 < \cdots < \beta_{n-1} <
\lam$.  Define $f \in {^\kappa \kappa} \cap
\calK$ as follows:
$$f(\nu) = \sup \{g(\nu_1\nek \nu_{k-1},
\nu) \mid \nu_1 < \cdots < \nu_{n-1} <
\nu\}\ .$$
\hfill$\bigsquare$

\proclaim Claim 2.  The set
$$A=\{\beta < i(\kappa) \mid \exists f
\in {^\kappa \kappa} \cap \calK, \quad i(f)
(\lam) = \beta\}$$
is $<\kappa$-closed.

\pr First suppose that the embedding $i$ is
definable in $\calK$.  Let $\{\beta_\nu\mid
\nu < \rho\}$ be a subset of $A$ for some
$\rho < \kappa$.  Denote by $\beta$ the
$\sup\{\beta_\nu\mid \nu < \rho\}$.  Since
there is no measurable cardinals of $\calK$
between $\kappa$ and $2^\kappa$, there will
be a set  $B \in \calK$, $|B| \le \kappa$
covering $\{\beta_\nu\mid \nu < \rho\}$, by
Mitchell [Mi2].  Notice that since $\kappa^+
= (\kappa^+)^\calK$, $|B|^\calK \le
\kappa$.  Using the regularity of $\kappa$,
we can find $B^* \subseteq B \cap A$, $B^*
\in \calK$ of cardinality $< \kappa$
cofinal in $\beta$.

Now it is obvious that for a function $f
\in {^\kappa \kappa} \cap \calK$,\enskip $i(f) (\lam)
= \beta$.  So $\beta \in A$.

Let us now deal with the general case, i.e.
$i$ is not necessarily definable in
$\calK$.  We replace $i$ by a definable in
$\calK$ iteration $i^*$.  Proceed as
follows.  Iterate $\calK$ using every
measure over $\kappa$ as well as the new
ones appearing in the process
$(2^\kappa)^+$--many times.  Let $i^*$ be
such an iteration.  Obviously, there is an
embedding $k: i(\calK) \to i^* (\calK)$ so
that the following diagram is commutative:
\bigskip
\input pictex
\font\thinlinefont=cmr5
\def\calK{{\cal K}}
$$\beginpicture
\setcoordinatesystem units < 0.500cm, 0.500cm>
%\unitlength= 1.000cm
\linethickness=1pt
\setplotsymbol 
%({\makebox(0,0)[l]{\tencirc\symbol{'160}}})
%\setshadesymbol 
({\thinlinefont .})
\setlinear
%
% Fig POLYLINE object
%
\linethickness= 0.500pt
\setplotsymbol ({\thinlinefont .})
\plot  6.953 19.082 10.287 20.987 /
%
% arrow head
%
\plot 10.098 20.806 10.287 20.987 10.035 20.916 /
%
%
% Fig POLYLINE object
%
\linethickness= 0.500pt
\setplotsymbol ({\thinlinefont .})
\plot  6.953 18.447 10.128 17.177 /
%
% arrow head
%
\plot  9.869 17.212 10.128 17.177  9.916 17.330 /
%
%
% Fig POLYLINE object
%
\linethickness= 0.500pt
\setplotsymbol ({\thinlinefont .})
\putrule from 11.239 20.510 to 11.239 17.812
%
% arrow head
%
\plot 11.176 18.066 11.239 17.812 11.303 18.066 /
%
%
% Fig TEXT object
%
\put {$i^*$} [lB] at  8.382 17.018
%
% Fig TEXT object
%
\put {$i$} [lB] at  8.223 20.510
%
% Fig TEXT object
%
\put {$\calK$} [lB] at  6.001 18.764
%
% Fig TEXT object
%
\put {$k$} [lB] at 11.874 19.082
%
% Fig TEXT object
%
\put {$i(\calK)$} [lB] at 10.763 20.987
%
% Fig TEXT object
%
\put {$i^*(\calK)$} [lB] at 10.604 17.177
%
% Fig TEXT object
%
\put { .} [lB] at 12.668 17.177
\linethickness=0pt
\putrectangle corners at  6.001 21.368 and 12.668 17.018
\endpicture$$
\noindent
Notice that $k \rhookup \kappa^+ = id$.
So, $k(A) = \{\beta< i^* (\kappa) |\enskip
\exists f \in {^\kappa \kappa} \cap \calK$,
$i^* (f) \bigl(k(\lam)\bigr) = \beta\} =
k^{''} (A)$.

Now, the argument used above for $i$, $A$ works
for $i^*$, $k(A)$.\hfill$\bigsquare$

Let us now change the cofinality of
$\kappa$ to $\omega$ by the Prikry forcing
with $U$.  Let $\langle \kappa_n \mid n <
\omega\rangle$ be the Prikry sequence and
$\langle \lam_n\mid n < \omega\rangle$ the
sequence corresponding to $\lam$, i.e.
$\langle f(\kappa_n) \mid n <
\omega\rangle$ for $f \in {^\kappa\kappa}$
representing $\lam$ in the ultrapower by
$U$.

Choose an elementary submodel $M$, of large
enough portion of the universe $|M| <
\kappa$, $^\omega M \subseteq M$ containing
all relevant information.  Let $\alp = \sup
\bigl(M \cap i(\kappa)\bigr)$ and $\alp_n =
\sup (M \cap \kappa_{n+1})$ for $n <
\omega$.  Notice that $\cf\alp <\kap$, since
$|M|<\kap$ and $\cf\alp >\ome$, since ${}\!^\ome
M\subseteq M$.  

Since $A = \{ \beta < i(\kappa)\mid \exists f
\in {^\kappa \kappa} \cap \calK$, $i(f)
(\lam) = \beta\}$ is $<\kappa$-closed and
unbounded in $i(\kappa)$ in $V$, it contains its
limit points of cofinality $\del,\ome <\del<\kap$,
in $V[\langle \kappa_n \mid n <
\omega\rangle ]$. Hence, $\alp \in A$,
since $A \in M$ and $\ome <\cf\alp <\kap$. 
Let $f_\alp \in {^\kappa \kappa} \cap \calK$ be such
that $i(f_\alp) (\lam) = \alp$.

\proclaim Claim 3.  For all but finitely
many $n$'s, $f_\alp (\lam_n) = \alp_n$.

\pr Let $g_\nu$  be a function in $V$  such that
$\nu=i(g_\nu)(\kap)$.  Then it is enough to show
that $g_\alp(\kap_n)=\alp_n$ for almost all $n$,
for then 
$$\eqalign{i(f_\alp)(\lam)=\alp&\Longleftrightarrow
\{ \nu :f_\alp(g_\lam (\nu))=g_\alp(\nu)\}\in U\cr
&\Longleftrightarrow f_\alp(g_\lam(\kap_n))=g_\alp
(\kap_n)\ \hbox{for sufficiently large}\ n\cr
&\Longleftrightarrow f_\alp(\lam_n)=\alp_n\
\hbox{for sufficiently large}\ n\ .\cr}$$

Now we show that $g_\alp(\kap_n)=\alp_n$  for
almost every $n$.  Let $g\in
M\cap\prodl_{n<\ome}\kap_{n+1}$.  By standard
argument on Prikry forcing used inside $M$,  we
can find $\bet\in M\cap\alp$  such that $g_\bet
(\kap_n)>g(n)$  for all but finitely many $n$'s.
So $\alp_n\le g_\alp(\kap_n)$.  
Now $\ome <\cf (\alp)<\kap$  in $V$  since $\ome
<\cf (\alp)<\kap$  in $V[(\kap_n:n\in\ome)]$.
Choose an increasing sequence $(\xi_\nu:\nu
<\rho)\in V$  which is cofinal in $\alp$,  for some
$\rho <\kap$.  Then
$$\left\{ \eta:g_\alp(\eta)=\sup\limits_{\nu
<\rho}g_{\xi_\nu}(\eta)\right\}\in U\ .$$
If $\alp_n<g_\alp(\kap_n)$  for finitely many
$n$  then there is $\nu <\rho$ such that
$\alp_n<g_{\xi_\nu}(\kap_n)$  for infinitely
many $n$  (using the fact that $\cf
(\alp)>\ome)$.  This is impossible, since if we
pick $\xi\in M$  with $\xi_\nu <\xi <\alp$ then  
$g_{\xi_\nu}(\kap_n)<g_\xi(\kap_n)<\alp_n$  for
almost all $n$.\hfill$\bigsquare$

Let us now use the club $E \in V$
consisting of regular in $\calK$ cardinals.
For every $n < \omega$, $E \cap
\kappa_{n+1}$ will be such a club in
$\kappa_{n+1}$.  Hence the same is true in $M$.  So,
$\alp_n = M \cap \kappa_{n+1} \in E$ for every
$n<\omega$.  We  argue that this is impossible. 

Briefly we apply Mitchell's analysis of $M$, see [Mi2]. 
It implies that for all but finitely many $n$'s,
$\alp_n$ is a limit of indiscernibles.
Since $\kappa_{n+1}$ is indiscernible for
$\kappa$, there will a club subset of
$\alp_n \cap E$ consisting of
indiscernibles for at least $\kappa$.  Select
a sequence $\langle c_n \mid n <
\omega\rangle$ of such indiscernibles so
that $\lam_n < c_n < \alp_n$.  Then
$\langle c_n \mid n < \omega\rangle \in M$
and by Claim 3, $f_\alp (\lam_n) > c_n$ $(n
< \omega)$,  which is impossible.
A contradiction.

Let us provide more details. 

We apply the technique of the second section of
[Mi2]. We will use $M$  as the set $N$  of that
paper, i.e. we are covering the set $M$.

Write $C(\alp)=\bigcup_\bet \calC(\alp,\bet)$
and $C=\bigcup_\alp C(\alp)$.  It follows
from what is given there that there are $h^M\in\calK$
and $\xi <\kap$  such that: 

\item{(1)(a)} $\forall\gam \in M\cap\kap\ \gam\in
h^M{}\!\tagg(\xi\cup(C\cap\gam +1))$. 
\item{(b)} If for $c\in C$  we write $\tau
(c)$  for the largest ordinal $\tau$ such that 
$c\in\calC(\tau)$  then $\tau (c)$  always
exists and is in $h^M{}\!\tagg (\xi\cup(C\cap c))$.

Since $\alp_n$ is regular in $\calK$,  it follows that
$C$ is unbounded in $\alp_n$ and $\kap_n\in C$ for
sufficiently large $n$.  
Since $(\kap_n:n\in\kap)$ is a Prikry sequence
it follows that $\tau (\kap_n)=\kap$ for
sufficiently large $n$.  Also it follows that
for $c$  in a closed unbounded subset of
$\alp_n \cap M$  we have $\tau(c)\not <\alp_n$.
Since $\tau(c)\in M$  it follows that
$\tau(c)>\alp_n$.   

An additional fact from [Mi2] that we are using:
\item{(2)} For every $f\in\calK$ there is $\rho <\kap$ 
such that if $c>\rho$  then $(f\tagg c)\cap\tau (c)
\subset c$. 

However we have we have $f_\alp(\lam_n)=\alp_n$  for all
sufficiently large $n$,  and this is impossible
because for all sufficiently large $n$  there are
$c_n$  such that $\lam_n<c_n<\alp_n$  and $\tau
(c_n)>\alp$.\hfill$\bigsquare$

We do not know if it is possible to improve the
above to $o(\kap)=\kap +1$. 

\subheading{Question} Is $o(\kap)=\kap$  enough
for a model with a measurable containing a club
of former regular cardinals?

\references{50}

\ref{[Av-Sh]} U. Avraham and S. Shelah, Forcing
closed unbounded sets, J. Sym. Logic 48 (1983),
643-648.
\smallskip
\ref{[Gi1]} M. Gitik, Changing cofinalities
and the nonstationary ideal, {\it Israel Journal
of Math.} 56(3) (1986),  280-314.
\smallskip
\ref{[Gi2]} M. Gitik, On measurable
cardinals violating GCH, {\it Ann. of Pure and
Appl. Logic} 63 (1993), 227-240.
\smallskip
\ref{[Gi3]} M. Gitik, On the Mitchell and
Rudin-Kiesler orderings of ultrafilters, Ann. of
Pure and appl. Logic 39 (1988), 175-197.
\smallskip
\ref{[Gi-Sh]} M. Gitik and S. Shelah, On Certain
Indestructibility of Strong Cardinals and a
Question of Hajnal, Arch. Math. Logic 28, (1989)
35-42. 
\smallskip
\ref{[Mi1]} W. Mitchell, How weak is a
closed unbounded ultrafilter?  In Logic
Coll. 80, North Holland 1982, pp. 209-231.
\smallskip
\ref{[Mi2]} W. Mitchell, Applications of
the Core Model for sequences of measures,
{\it Trans. Amer. Math. Soc.} 299 (1987), 41-58.
\smallskip
\ref{[Mi3]} W. Mitchell,  The Core Model for
sequences of measures I, Math. Proc. Cambridge
Phil. Soc. 95 (1984), 229-260.

\end
\input pictex
\font\thinlinefont=cmr5
\def\calK{{\cal K}}
$$\beginpicture
\setcoordinatesystem units < 0.500cm, 0.500cm>
%\unitlength= 1.000cm
\linethickness=1pt
\setplotsymbol 
%({\makebox(0,0)[l]{\tencirc\symbol{'160}}})
%\setshadesymbol 
({\thinlinefont .})
\setlinear
%
% Fig POLYLINE object
%
\linethickness= 0.500pt
\setplotsymbol ({\thinlinefont .})
\plot  6.953 19.082 10.287 20.987 /
%
% arrow head
%
\plot 10.098 20.806 10.287 20.987 10.035 20.916 /
%
%
% Fig POLYLINE object
%
\linethickness= 0.500pt
\setplotsymbol ({\thinlinefont .})
\plot  6.953 18.447 10.128 17.177 /
%
% arrow head
%
\plot  9.869 17.212 10.128 17.177  9.916 17.330 /
%
%
% Fig POLYLINE object
%
\linethickness= 0.500pt
\setplotsymbol ({\thinlinefont .})
\putrule from 11.239 20.510 to 11.239 17.812
%
% arrow head
%
\plot 11.176 18.066 11.239 17.812 11.303 18.066 /
%
%
% Fig TEXT object
%
\put {$i^*$} [lB] at  8.382 17.018
%
% Fig TEXT object
%
\put {$i$} [lB] at  8.223 20.510
%
% Fig TEXT object
%
\put {$\calK$} [lB] at  6.001 18.764
%
% Fig TEXT object
%
\put {$k$} [lB] at 11.874 19.082
%
% Fig TEXT object
%
\put {$i(\calK)$} [lB] at 10.763 20.987
%
% Fig TEXT object
%
\put {$i^*(\calK)$} [lB] at 10.604 17.177
%
% Fig TEXT object
%
\put { .} [lB] at 12.668 17.177
\linethickness=0pt
\putrectangle corners at  6.001 21.368 and 12.668 17.018
\endpicture$$